\documentclass[11pt]{amsart}
\usepackage{amssymb}
\usepackage{graphicx}
\usepackage{amscd}
\usepackage{url}
\usepackage[utf8]{inputenc}
\usepackage[T1]{fontenc}

\theoremstyle{plain}
\newtheorem{theorem}{Theorem}
\newtheorem*{acknowledgement}{Acknowledgement}

\newtheorem{lemma}[theorem]{Lemma}

\theoremstyle{definition}

\newtheorem*{example}{Example}
\newtheorem*{algorithm}{Algorithm}

\DeclareMathOperator{\ord}{ord}

\title[Puiseux expansion]{An efficient method of finding a Puiseux expansion of a parametric singularity}
\author{Maciej Borodzik}
\address{Institute of Mathematics, University of Warsaw, ul. Banacha 2,
02-097 Warsaw, Poland}
\email{mcboro@mimuw.edu.pl}
\date{8 Dec 2009}
\thanks{Supported by Polish  MNiSz Grant N N201 397937. The author is also supported by Foundation for Polish Science (FNP)}
\begin{document}
\begin{abstract}
We present a new and very efficient method of determining a topological type of a cuspidal plane curve singularity parametrised
locally by $(x(t),y(t))$. The approach
relies on defining inductively a sequence of functions $P_k$ which are polynomials in $x$, $y$ and derivatives of $x$ and $y$.
\end{abstract}
\maketitle

Let us be given a germ of a parametric curve in $\mathbb{C}^2$ given by
\begin{equation}\label{eq:param}
\begin{split}
x(t)&=a_0t^p+a_1t^{p+1}+\ldots,\\
y(t)&=b_0t^q+b_1t^{q+1}+\ldots.
\end{split}
\end{equation}
Here $a_0b_0\neq 0$, $q\ge p>1$ and both series on the right hand side are assumed to be convergent in some neighbourhood of $t=0$.
In order to determine the topological type of the singularity \eqref{eq:param} we need (see \cite{BK,EN})
to find its characteristic sequence
$(p;q_0,q_1,\ldots,q_n)$. To do this we write the Puiseux expansion of \eqref{eq:param} in the form
\begin{equation}\label{eq:puis}
y=c_qx^{q/p}+c_{q+1}x^{(q+1)/p}+\ldots.
\end{equation}
Then we define $q_0=\min\{s\ge 0\colon c_s\neq 0 \wedge p\nmid s\}$, and inductively $p_0=p$, $p_k=\gcd(p_{k-1},q_{k-1})$. 
If $p_k\neq 1$ we define $q_k=\min\{s\ge 0\colon c_s\neq 0\wedge p_k\nmid s\}$; if $p_k=1$ we stop the procedure. As $p_k>p_{k+1}$ this
procedure will definitely stop after a finite number of steps. So let $n$ be such that $p_{n+1}=1$, then
$(p;q_0,\ldots,q_n)$ is the characteristic sequence.

Therefore, in order to know the topological type of singularity we need to know the coefficients $c_k$, or, at any rate, 
to know which one
are zero and which are not. Now passing from \eqref{eq:param} to \eqref{eq:puis} is, in general, a very involved procedure, which requires
finding an explicit formula for $t=d_1x^{1/p}+d_2x^{2/p}+\ldots$. Moreover, the Puiseux coefficients $c_k$ are not very well defined, 
namely they are functions of some fractional powers of $a_0$. In particular if we study a deformation, in other words
if we allow the coefficients
of \eqref{eq:param} to vary with some parameter $s$, and $a_0(s)\to 0$ as $s\to s_0$, then $c_k(s)$ can easily get 
out of hand. Some of them can diverge to
infinity, some of them can converge but lose any topological meaning. We want to propose a remedy to both problems: we shall define some 
polynomials in $y$, $x$ and their derivatives over $t$, which carry the same data about the singularity as Puiseux expansion and such
that one can recover the Puiseux expansion from them in a very easy way. 

As a matter of fact, the definition of polynomials $P_k$ appears first in \cite{Bo}, but only after writing this article the author
realised that the polynomials $P_k$ may be of wider use and allow fast and efficient computation of a characteristic sequence.

We begin with putting 
\[P_0(t)=y\text{ and }r_0=\ord_{t=0}P_0(t)=q\] 

Like in \cite[Proof of Lemma 3.2]{BZ}, let us divide \eqref{eq:puis} by $x^{q/p}$ and 
differentiate both sides with respect to $t$ (it is convenient to write here
$\dot x=\frac{dx}{dt}$ and so on). We get
\[\frac{\dot y x-\frac{q}{p}\dot x y}{x^{q/p+1}}=\frac{1}{p}c_{q+1}\dot x x^{1/p-1}+\frac{2}{p}c_{q+2}\dot x x^{2/p-1}+\ldots.\]
So if
\begin{equation}\label{eq:p1}
P_1=\dot yx-\frac{q}{p}\dot x y
\end{equation}
we will have
\[P_1=\frac{1}{p}c_{q+1}x^{(q+1)/p}\dot x+\frac{2}{p}c_{q+2}x^{(q+2)/p}\dot x+\dots.\]
Let 
\begin{equation}\label{eq:r1}
r_1=\ord_{t=0}P_1-(p-1).
\end{equation}
Then we have an obvious lemma.
\begin{lemma}
The coefficients $c_{r_0+1},\ldots,c_{r_1-1}$ all vanish.  Moreover the first not vanishing coefficient of $P_1$, i.e. the one at $t^{r_1+(p-1)}$,
is equal to $(r_1-r_0)c_{r_1}\cdot a_0^{r_1/p+1}$.
\end{lemma}
Now we define inductively the polynomials $P_k$ and the numbers $r_k$. The induction procedure is as follows. Assume that
we have
\begin{equation}\label{eq:pk}
P_k=\frac{r_k-r_{k-1}}{p}\cdot \frac{r_k-r_{k-2}}{p}\cdot\ldots\cdot \frac{r_k-r_0}{p}c_{r_k}\dot x^{2k-1}x^{r_k/p}+\ldots
\end{equation}
with
\begin{equation}\label{eq:rk}
r_k=\ord_{t=0}P_k-(2k-1)(p-1).
\end{equation}
In \eqref{eq:pk} dots denote terms like $(s-r_{k-1})\ldots (s-r_0)p^{-k}c_s\dot x^{2k-1}x^{s/p}$, which are too long
to be written in a nice looking formula.

Now we divide both sides of \eqref{eq:pk} by ${\dot x}^{2k-1}x^{r_k/p}$, differentiate both sides with respect to $t$
and then multiply them again by $\dot x^{2k}x^{r_k/p+1}$. We get
\[
P_{k+1}=\frac{r_{k+1}-r_k}{p}\frac{r_{k+1}-r_{k-1}}{p}\cdot\ldots\cdot \frac{r_{k+1}-r_0}{p}c_{r_{k+1}}
\dot x^{2k+1}x^{ r_{k+1}/p}+\ldots,
\]
where
\begin{equation}\label{eq:rk1}
r_{k+1}=\ord_{t=0}P_{k+1}-(2k+1)(p-1)
\end{equation}
and
\begin{equation}\label{eq:induct}
P_{k+1}=x\dot x\frac{d}{dt}P_k-\left(\frac{r_k}{p}{\dot x}^2+(2k-1)\ddot x x\right)P_k.
\end{equation}
Similarly we have a lemma the proof of which is an easy computation.
\begin{lemma}\label{lem:2}
The Puiseux coefficients $c_{r_{k-1}+1},\ldots,c_{r_k-1}$ all vanish. The coefficient at $t^{r_k+(2k-1)(p-1)}$ of $P_k$, i.e. the
first non-vanishing one,
is equal to
\[(r_k-r_{k-1})(r_k-r_{k-2})\ldots (r_k-r_0)p^{k-1}c_{r_k}a_0^{r_k/p+(2k-1)}.\]
\end{lemma}
Thus, using the above method we can easily find the Puiseux expansion. The computations become simpler if there are many
Puiseux terms that vanish.
The algorithm, in a form easy to use, can be presented as follows.
\begin{algorithm}
Given a parametric presentation as in \eqref{eq:param} of a cuspidal singularity with $a_0b_0\neq 0$.
\begin{itemize}
\item[Step 1.] Set $p_0=p$ and  $r_0=q$. Let $p_1=\gcd(p_0,q_0)$. If $p_1=1$ then the singularity is quasi-homogeneous and we stop.
\item[Step 2.] Compute $P_1$ as in \eqref{eq:p1} and $r_1$ from \eqref{eq:r1}. Let $p_2=\gcd(p_1,r_1)$.
\item[Step 3.] If $p_2=1$ then we stop. Otherwise put $k=1$.
\item[Step 4.] Compute $P_{k+1}$ and $r_{k+1}$ from \eqref{eq:induct} and \eqref{eq:rk1} respectively.
\item[Step 5.] Let $p_{k+1}=\gcd(p_k,r_k)$ if $p_{k+1}=1$ we stop, otherwise we increase $k$ and go to Step 4.
\end{itemize} 
As an output we get a sequence of positive integers $(r_0,r_1,r_2,\dots,r_n)$, such that the Puiseux expansion is
\[y=c_{r_0}x^{r_0/p}+c_{r_1}x^{r_1/p}+\dots+c_{r_k}x^{r_k/p}+\textrm{inessential Puiseux terms.}\]
Here $c_{r_0},c_{r_1},\ldots,c_{r_k}$ are all non-zero.
The procedure will eventually stop unless the parametrisation \eqref{eq:param} is not one to one.
\end{algorithm}
Of course, we can recover not only the characteristic sequence, but also the exact values of non-vanishing
Puiseux coefficients. It is enought to look at the first non-vanishing coefficients of functions $P_{k+1}$
and use second part of Lemma~\ref{lem:2}. 
\begin{example}
Consider a singularity parametrised by
\begin{align*}
x(t)&=t^{12}+t^{13}+\frac{37}{28}t^{14}\\
y(t)&=t^{18}+\frac{3}{2}t^{19}+\frac{33}{14}t^{20}+\frac{13}{14}t^{21}+\frac{675}{1568}t^{22}-\frac{675}{3136}t^{23}.
\end{align*}
We write
\[y=x^{18/12}+c_{19}x^{19/12}+\ldots.\]
But it is easy to compute that
\[P_1=-\frac{2025}{10976}t^{35}-\frac{24975}{43904}t^{36}.\]
So $r_1=35-11=24$. Thus $c_{19}=c_{20}=c_{21}=c_{22}=c_{23}=0$ and $c_{24}\neq 0$. Then
\[P_2=\frac{2500875}{76832}t^{59}+\ldots\]
so $r_2=59-3\cdot 11=26$, so $c_{25}=0$ and $c_{26}\neq 0$. 
Then a simple computation shows that $P_3\sim t^{82}$, so $c_{27}\neq 0$ and
the Puiseux expansion is $y=x^{3/2}+c_{24}x^2+c_{26}x^{13/6}+c_{27}x^{27/12}+\ldots$, 
so the characteristic sequence is $(12;18,26,27)$ ($c_{24}$ is inessential).
\end{example}

We can use the polynomials $P_k$ in studying deformations. In fact, they behave very well if $x(t)$ and $y(t)$
vary smoothly. This approach is pursued in \cite{Bo} with quite a success in some partial cases.

\begin{acknowledgement}
The author wishes to express his thanks for H.~Żołądek, A.~Płoski and P.~Cassou-Nogu\`es 
for fruitful discussions and to Christian Gorzel
for his many comments to the first version of the paper. 
\end{acknowledgement}

\end{document}